\documentclass[12pt]{amsart}
\usepackage{amsmath,amssymb}

\newtheorem{theorem}{Theorem}[section]
\newtheorem{lemma}[theorem]{Lemma}
\newtheorem{claim}[theorem]{Claim}

\newtheorem{question}[theorem]{Question}

\theoremstyle{definition}
\newtheorem{definition}[theorem]{Definition}
\newtheorem{example}[theorem]{Example}

\theoremstyle{remark}
\newtheorem{remark}[theorem]{Remark}

\DeclareMathOperator{\V}{V} \DeclareMathOperator{\osc}{osc}

\DeclareMathOperator{\diam}{diam} \DeclareMathOperator{\Vol}{Vol}

\author[Z.L. Zhang]{Zhenlei Zhang}
\address{Department of Mathematics\\
Capital Normal University\\  Xisanhuan North Road 105\\ Beijing,
100048, P.R.China\\}

\email{zhleigo@yahoo.com.cn}


\begin{document}

\title{K\"ahler Ricci flow with vanished Futaki invariant}

\begin{abstract}
We study the convergence of the K\"ahler-Ricci flow on a compact K\"ahler manifold $(M,J)$ with positive first Chern class $c_1(M;J)$ and vanished Futaki invariant on $\pi c_1(M;J)$. As the application we establish a criterion for the stability of the K\"ahler-Ricci flow (with perturbed complex structure) around a K\"ahler-Einstein metric with positive scalar curvature, under certain local stable condition on the dimension of holomorphic vector fields. In particular this gives a stability theorem for the existence of K\"ahler-Einstein metrics on a K\"ahler manifold with possibly nontrivial holomorphic vector fields.
\end{abstract}

\today

\maketitle

\section{Introduction}

This is a continuation of \cite{Zhang}. We shall apply the results of \cite{Zhang} to study the convergence and stability of the K\"ahler-Ricci flow with positive first Chern class and vanished Futaki invariant. The convergence in this paper is always taken with respect to a fixed complex structure under the K\"ahler-Ricci flow. The results relate closely to Chen and Li's theorems in \cite{ChLi}.

To start our argument, recall that for any compact K\"ahler manifold $(M,g,J)$ with $c_1(M;J)>0$, if its K\"ahler form lies in $\pi c_1(M;J)$, then through the $\partial\bar{\partial}$-lemma, there exists a unique Ricci potential $u$, a real-valued function, in the sense of
\begin{equation}\label{RP}
R_{i\bar{j}}+\partial_i\partial_{\bar{j}}u=g_{i\bar{j}}
\end{equation}
with the normalization
\begin{equation}\label{RPN}
\frac{1}{\V}\int_Me^{-u}dv=1
\end{equation}
where $\V=\int_Mdv$ denotes the volume of the manifold.

Define the weighted Laplace operator $L=-\triangle+g^{i\bar{j}}\nabla_iu\nabla_{\bar{j}}$. It is a self-adjoint operator on the space $L^2(e^{-u}dv)$ whose smallest positive eigenvalue is $\geq 1$ (cf. \cite{Fu} or \cite{TiZh2}), with equality holds iff there exists a nontrivial holomorphic vector field. What we are interested in is the ``second" eigenvalue of $L$. Namely, denote by $\lambda(g)$ the first positive eigenvalue
of $L$ if $(M,J)$ has no nonzero holomorphic vector fields;
otherwise, $\lambda(g)$ denotes the second positive eigenvalue of
$L$.

Let $(M,g_0,J)$ be a compact K\"ahler manifold with $c_1(M;J)>0$.
Suppose its K\"ahler form belongs to $\pi c_1(M)$ and the Futaki invariant vanishes on $\pi c_1(M;J)$. Consider the K\"ahler-Ricci flow with initial metric $g(0)=g_0$:
\begin{equation}\label{KRF}
\frac{\partial}{\partial t}g_{i\bar{j}}=-R_{i\bar{j}}+g_{i\bar{j}},
\end{equation}
where $R_{i\bar{j}}$ denotes the Ricci curvature of $g$. It is obvious that the K\"ahler class and the volume are preserved
by the K\"ahler-Ricci flow. By Cao \cite{Ca}, the solution exists for all time $t\in[0,\infty)$. In \cite{Zhang}, the author showed how to control the convergence of the K\"ahler-Ricci flow when $\lambda(g(t))$ is uniformly lower bounded away from $1$. In this paper, we aim to bound $\lambda(g(t))$ from below in terms of the geometry of the initial metric $g(0)$ and thus obtain a global convergence theorem of $g(t)$. Our main theorem reads

\begin{theorem}\label{main theorem}
Given positive constants $\delta$, $\Lambda$ and $D$, there exists
$\epsilon$ such that for any K\"ahler manifold $(M,g_0,J)$ whose
K\"ahler form lies in $\pi c_1(M;J)$, if the Futaki invariant vanishes
on $\pi c_1(M;J)$ and
\begin{equation}\label{Z1}
bisc.(M,g_0)\geq-\Lambda,\hspace{0.2cm}\diam(M,g_0)\leq D,\hspace{0.2cm}\lambda(g_0)\geq 1+\delta,\hspace{0.2cm}\|\nabla u_0\|_{C^0}\leq\epsilon,
\end{equation}
then the K\"ahler-Ricci flow starting from $(g_0,J)$ will converge
exponentially fast to a K\"ahler-Einstein metric $(g_\infty,J)$.
\end{theorem}

Here, $bisc.(M,g_0)$ denotes the bisectional curvature of
the K\"ahler manifold $(M,g_0)$, while $u_0$ denotes the K\"ahler
potential of $(M,g_0)$ determined by (\ref{RP}) and (\ref{RPN}).

\begin{theorem}\label{t2}
Given positive constants $n,\delta,D$ and $\Lambda$, there exists
$\epsilon$ such that for any K\"ahler manifold $(M,g_0,J)$ of dimension $n$ whose
K\"ahler class is $\pi c_1(M;J)$, if the Futaki invariant vanishes
on $\pi c_1(M;J)$ and
\begin{equation}\label{Z2}
|Rm(g_0)|\leq\Lambda,\hspace{0.2cm}\lambda(g_0)\geq 1+\delta,\hspace{0.2cm}\frac{1}{\V}\int_M(s-n)^2dv\leq\epsilon,
\end{equation}
where $Rm$ and $s$ denote the Riemannian curvature tensor and scalar curvature respectively, then the K\"ahler-Ricci flow starting from $(g_0,J)$ will converge exponentially fast to a K\"ahler-Einstein metric $(g_\infty,J)$.
\end{theorem}

An immediate consequence of above theorem is the following criterion theorem for the stability of the K\"ahler-Ricci flow around a K\"ahler-Einstein metric. Notice that for any metric $g$ which is $C^2$ close to a
K\"ahler-Einstein metric $g_{KE}$, its geometry is uniformly bounded
and its traceless Ricci curvature is sufficiently small.

\begin{theorem}\label{t3}
Let $(M,g_{KE},J_{KE})$ be a K\"ahler-Einstein manifold with
positive scalar curvature. Then for any $\delta>0$, there exists
$\epsilon>0$ such that for any $\epsilon$-deformed K\"ahler manifold $(M,g,J)$ in the sense that
\begin{equation}\label{Z3}
\|g-g_{KE}\|_{C^2(g_{KE})}+\|J-J_{KE}\|_{C^2(g_{KE})}\leq\epsilon,\hspace{0.3cm}\lambda(g)\geq
1+\delta,
\end{equation}
if the K\"ahler form of $g$ belongs to $\pi c_1(M;J)$ and the Futaki invariant
vanishes on $\pi c_1(M;J)$, then the K\"ahler-Ricci flow starting from $(g,J)$ will converge
exponentially fast to a K\"ahler-Einstein metric $(g_\infty,J)$.
\end{theorem}

\begin{remark}
If the perturbed complex structure $J$ in the theorem is stable in
the sense of Condition (B) in \cite{PhSt} or Definition 4.14 in
\cite{ChLiWa}, then the eigenvalue $\lambda(g)$ is uniformly lower bounded
away from $1$ when $\epsilon$ is small enough. The Theorem 3 of \cite{PhSt} gives a proof of this fact when $(g,J)$ is $C^\infty$ close to $(g_{KE},J_{KE})$. When $(g,J)$ is only $C^2$ close to $(g_{KE},J_{KE})$, we can run the Ricci flow for a short time to get the required $C^\infty$ approximation. Based on this discussion, we get another proof of the stability theorem, namely Theorem 1.6, of Chen and Li
\cite{ChLi}; the Theorem \ref{t2} can also be viewed as a generalization of Theorem 1.10 of \cite{ChLi}.
\end{remark}

\begin{remark}
Donaldson \cite{Do} showed the existence of a K\"ahler-Einstein metric
$(g_{KE},J_{KE})$ on the Mukai-Umemura manifold; however, according
to Tian's famous counterexample \cite{Ti2}, the generic deformation
of the Mukai-Umemura manifold, which has no nontrivial holomorphic
vector fields, does not admit any K\"ahler-Einstein metric. In view of
above theorem, it means that the first positive eigenvalue $\lambda(g)$
of $(g,J)$, a small deformation of $(g_{KE},J_{KE})$, is arbitrarily
close to $1$ and the gradient of its eigenfunctions give rise to the holomorphic
vector fields of $J_{KE}$ up to taking limits.
\end{remark}

According to Tian and Donaldson's counterexample, the jumping phenomena of the dimension of holomorphic vector fields is essential for the stability of the existence of K\"ahler-Einstein metrics on a nearby complex structure. In the following, we shall show that these two stability condition, namely the dimension of holomorphic vector fields and the existence of K\"ahler-Einstein metrics, are somehow equivalent, at least under the condition of vanished Futaki invariant. In other words, if the dimension of holomorphic vector fields is ``stable", then the existence of K\"ahler-Einstein metrics should also be ``stable".

To make the notion clear, we introduce the following definition. Denote by $h^0(J)$ the space of holomorphic vector fields with respect to the complex structure $J$.

\begin{definition}
A complex structure $J$ on a complex manifold $M$ is called
``relative stable", if $\dim h^0$,
defined on the space of complex structures on $M$, is a local
minimizer at $J$.
\end{definition}

In particular, $J$ will always be ``relative stable" if $(M,J)$ has
no nontrivial holomorphic vector fields. Moreover, because $\dim h^0$ is lower semi-continuous on the space of complex structures, the condition ``relative stable" is an open condition and $J$ being ``relative stable" implies that $\dim h^0$ is actually constant on a small neighborhood of $J$.

\begin{example}
As showed in Lemma 1.2 of \cite{Ti0}, any complex surfaces of the
form $\mathbb{C}P^2\sharp n\overline{\mathbb{C}P^2}$ ($5\leq n\leq
8$) with positive first Chern class has no nontrivial holomorphic
vector field and hence is ``relative stable".
\end{example}

Let $(M,g,J)$ be a K\"ahler manifold with $c_1(M;J)>0$. It is an easy check that if the complex structure $J$ is ``relative stable", then
for any small deformation, say $(g^{'},J^{'})$, of $(g,J)$ in the $C^2$ topoloty, the eigenvalue
$\lambda(g^{'})$ admits a uniform lower bound away from $1$. Hence we can
prove the following theorem.

\begin{theorem}\label{t4}
Let $(M,g_{KE},J_{KE})$ be a K\"ahler-Einstein manifold with
positive scalar curvature such that $J_{KE}$ is ``relative stable".
There exists $\epsilon>0$ such that for any $\epsilon$-deformed K\"ahler manifold
$(M,g,J)$ in the sense that
\begin{equation}\label{Z3}
\|g-g_{KE}\|_{C^2(g_{KE})}+\|J-J_{KE}\|_{C^2(g_{KE})}\leq\epsilon,
\end{equation}
if the K\"ahler form of $g$ belongs to $\pi c_1(M;J)$ and Futaki invariant
vanishes on $\pi c_1(M;J)$, then the K\"ahler-Ricci flow starting from $(g,J)$ will converge
exponentially fast to a K\"ahler-Einstein metric $(g_\infty,J)$.
\end{theorem}

\begin{remark}
The special case when $(M,J_{KE})$ has no holomorphic vector fields has already been proved independently by Chen and Li in \cite{ChLi}.
\end{remark}

\begin{remark}
The relationship between the dimension of holomorphic vector fields and the convergence of a K\"ahler-Ricci flow was first studied by Phong and Sturm \cite{PhSt}. They introduced a different stability condition there, on the dimension of holomorphic $T^{1,0}$ vector fields under diffeomorphic actions up to taking Cheeger-Gromov-Hausdorff limits, to control the convergence of the K\"ahler-Ricci flow.
\end{remark}

\begin{remark}
Combining with the convergence theorem of Perelman, or Tian and Zhu \cite{TiZh2}, for the K\"ahler-Ricci flow, once known the existence of K\"ahler-Einstein metric in above theorem, the K\"ahler-Ricci flow will converge automatically (up to holomorphic transformations) for any initial metric $g$ whenever its its K\"ahler form lies in $\pi c_1(M;J)$.
\end{remark}

Finally let us briefly discuss the known convergence results previously for K\"ahler-Ricci flow on a compact manifold with positive first Chern class.  Suppose the existence of a K\"ahler-Einstein metric, Chen and Tian \cite{ChTi1,
ChTi2} proved the convergence to this metric when the bisectional curvature is positive; later, Perelman
showed the convergence without any curvature condition, and this was later extended to shrinking
K\"ahler-Ricci solitons by Tian and Zhu \cite{TiZh2}. Without assuming the existence of K\"ahler-Einstein metrics, in
\cite{PSSW1, PSSW2, PhSt}, Phong, Song, Sturm and Weinkove studied the K\"ahler-Ricci flow with two stability conditions-lower bounded Mabuchi
K-energy and lower bound of first eigenvalue of $\bar{\partial}$ operator on $T^{1,0}$ vector fields-and proved the convergence to K\"ahler-Einstein
metrics; the result was later generalized to shrinking K\"ahler-Ricci soliton cases by them \cite{PSSW3}. In \cite{ChLiWa}, Chen, Li and Wang considered the K\"ahler-Ricci flow with different stability conditions-lower bounded energy $E_1$ which was introduced in \cite{ChTi1} and pre-stability of the complex structure-and proved the convergence result when $E_1$ is almost minimizing; in a subsequent paper \cite{ChLi}, Chen and Li improved the condition on energy $E_1$, to the vanished Futaki invariant and small Calabi functional, to derive the convergence; stability theorems are also established for K\"ahler-Ricci flow in \cite{ChLi} when the complex structure is pre-stable and the Futaki invariant vanishes. Several other stability theorems for K\"ahler-Ricci flow around a K\"ahler-Einstein metric or generally shrinking K\"ahler-Ricci soliton have
also been established independently through different methods by Tian and Zhu \cite{TiZh3}, Zhu \cite{Zhu}, Sun and
Wang \cite{SuWa}, Zheng \cite{Zheng}; in these works, the convergence are always taken up to holomorphic transformations. Besides these, we refer to \cite{CaZh, Ch, ChLi0, ChWa1, ChWa2, ChWa3, MuSz, PhSeSt, RuZhZh, Se, SeTi, Sz, To, Zhu1} for other works on the convergence of the K\"ahler-Ricci flow with positive first Chern class.

The rest of the paper is organized as follows: In \S 2, we prove the main convergence Theorem
\ref{main theorem}; In \S 3, we give a short proof of Theorem \ref{t2}; In \S 4, we give some remarks and questions on the stability of existence of K\"ahler-Einstein metrics on a compact K\"ahler manifold with positive first Chern class; In the Appendix \S 5, we prove a gradient estimate to the eigenfunctions of the weighted Laplace $L$.

\medskip

\noindent {\bf Acknowledgement:} The author would like to thank Professors J. Song and B. Weinkove for their suggestions on Tian and Donaldson's counterexample. The author also would like to thank Professors F.Q. Fang and Y.G. Zhang for their discussions in the course of writing the paper. The author also thanks V. Tosatti for his nice comments on Remark 4.1 and suggestions on the work of Stoppa and Sz\'ekelyhidi.

\section{Proof of Theorem \ref{main theorem}}

In this section, $(M,g_0,J)$ denotes a compact K\"ahler manifold of
complex dimension $n$ with $c_1(M;J)>0$. We assume the K\"ahler
class of $g_0$ is $\pi c_1(M;J)$. Suppose that the Futaki invariant
vanishes on $\pi c_1(M,J)$ and $g_0$ satisfies condition (\ref{Z1}):
\begin{equation}\label{C1}
bisc.(M,g_0)\geq-\Lambda,\hspace{0.2cm}\diam(M,g_0)\leq
D,\hspace{0.2cm}\lambda(g_0)\geq 1+\delta,\hspace{0.2cm}\|\nabla
u_0\|_{C^0}\leq\epsilon,
\end{equation}
for fixed positive constants
$\delta,\Lambda, D$ and suitably chosen $\epsilon$. Let
$g(t),t\in[0,\infty)$, be the K\"ahler-Ricci flow starting from
$g_0$. The aim is to show that $g(t)$ converges exponentially fast
to a K\"ahler-Einstein metric whenever $\epsilon$ is chosen small
enough.

Let $u(t)$ be associated Ricci potential of $g(t)$ determined by (\ref{RP}) and (\ref{RPN}) and let $a(t)$ be its average defined by \begin{equation}\label{average}
a=\frac{1}{\V}\int_Mue^{-u}dv
\end{equation}
at each time, where $\V$ denotes the volume of the K\"aher-Ricci flow. As in Introduction, let $\lambda(t)=\lambda(g(t))$ be the ``second" eigenvalue of $L=-\triangle+g^{i\bar{j}}\nabla_iu\nabla_{\bar{j}}$ at time $t$. We may write $\lambda(t)=1+\delta(t)$ for a family of constants $\delta(t)>0$.

\subsection{Preliminary lemmas}

Before the proof of Theorem \ref{main theorem},  let us first establish several lemmas.

Define a function
\begin{equation}\label{e20}
\delta^{'}=\delta^{'}(a,s)=\frac{a}{1+e^{s}+a
e^{s}},\hspace{0.3cm}\forall a>0,s>0.
\end{equation}
Applying the Proposition 3.1 in \cite{Zhang} to our K\"ahler-Ricci flow $g(t)$ gives

\begin{lemma}\label{l21}
The following estimate for Ricci potential holds at any time $t$
\begin{equation}\label{RPP}
\int_M|\nabla u|^2e^{-u}dv\geq(1+\delta^{'})\int_M(u-a)^2e^{-u}dv
\end{equation}
where $\delta^{'}=\delta^{'}(\delta,\osc(u))$.
\end{lemma}

Note that $\delta^{'}$ increases in $a$ and decrease in $s$. Thus, (\ref{RPP}) remains true for any $\delta^{'}(a,s)$ with $a\leq\delta$ and $s\geq\osc(u)$.

Under the K\"ahler-Ricci flow $g(t)$, introduce for any time $t\in[0,\infty)$,
\begin{equation}\label{e21}
Y=\frac{1}{\V}\int_M(u-a)^2e^{-u}dv.
\end{equation}
The second lemma shows the exponential decay of $Y$.

\begin{lemma}\label{l22}
Suppose that {\rm(\ref{RPP})} holds for some uniform $\delta^{'}$. If
\begin{equation}\label{e22}
\|u-a\|_{C^0}\leq\min(\frac{1}{4},\frac{\delta^{'}}{8})
\end{equation}
on some time interval $[0,T]$, then we have
\begin{equation}\label{e23}
Y(t)\leq e^{-\delta^{'} t}Y(0),\hspace{0.3cm}\forall t\in[0,T].
\end{equation}
\end{lemma}
\begin{proof}
By the computation in \cite{Zhang},
\begin{eqnarray}
\frac{d}{dt}Y\leq\big((-2+2\|u-a\|_{C^0})(1+\delta^{'})+(2+\|u-a\|_{C^0})\big) Y\leq-\delta^{'} Y,\nonumber
\end{eqnarray}
whenever (\ref{e22}) holds. The result then follows.
\end{proof}

Next we show the bound of $\|u-a\|_{C^0}$ in terms of $Y$. To this purpose, we introduce a definition of local volume non-collapsing.

\begin{definition}
Given $\kappa>0$, a compact Rieamnnian manifold is called $\kappa$-non-collapsed on scales $\leq\rho$ if any metric ball $B(r)$ of radius $r\leq\rho$ satisfies the bound
\begin{equation}
\frac{\Vol(B(r))}{\Vol(M)}\geq\kappa r^m
\end{equation}
on volume ratio, where $m$ is the (real) dimension of the manifold.
\end{definition}

In view of Bishop-Gromov relative volume comparison theorem, our initial manifold $(M,g_0)$ is always $\kappa_0$-non-collapsed on scales $\leq 1$ for some constant $\kappa_0$ depending only on $n,\Lambda$ and $D$.

\begin{lemma}\label{l23}
Under the K\"ahler-Ricci flow, if $g(t)$ is $\kappa$-non-collapsed on scales $\leq\rho$ for some $\rho\leq 1$ and
\begin{equation}\label{e23.5}
\|\nabla u\|_{C^0}\leq 1,\hspace{0.5cm}\|u-a\|_{C^0}\leq 2\rho,
\end{equation}
then,
\begin{equation}\label{e24}
\|u-a\|_{C^0}\leq K\cdot\kappa^{\frac{-1}{2n+2}}\cdot Y^{\frac{1}{2n+2}},
\end{equation}
for a universal constant $K$.
\end{lemma}
\begin{proof}
We adopt the argument as in Lemma 3 of \cite{PSSW1}. Let $x$ be a point such that $A=|u-a|(x)=\|u-a\|_{C^0}$. Then $|u-a|\geq\frac{A}{2}$ on the metric ball $B=B(x,\frac{A}{2})$ of radius $\frac{A}{2}\leq\rho$. Thus, by definition of $\kappa$-non-collapsed,
$$Y\geq\frac{1}{\V}\int_B(u-a)^2e^{-u}dv\geq e^{-\max (u)}(\frac{A}{2})^2\kappa(\frac{A}{2})^{2n}.$$
From this we obtain
$$A\leq 2e^{\frac{\max(u)}{2n+2}}\cdot\kappa^{\frac{-1}{2n+2}}\cdot Y^{\frac{1}{2n+2}},$$
Since $a\leq 0$ by Jensen inequality, see \cite[Remark 5.1]{Zhang}, we get
$$\max(u)\leq a+\|u-a\|_{C^0}\leq 2.$$ The result then follows.
\end{proof}

\subsection{A uniform $C^2$ estimate of K\"ahler potential}

This subsection is devoted to proving a uniform $C^2$ estimate of K\"ahler potential in terms of its $C^0$ bound under the K\"ahler-Ricci flow. The method is introduced by Yau in his resolution of Calabi conjecture \cite{Ya} for solving an elliptic Monge-Amp\`ere equation. The parabolic version is given by Cao \cite{Ca}. Here, we use the argument of \cite{ChLiWa}.

Let $\phi(\cdot,t)$ be the solution to the differential equation
\begin{equation}\label{e25}
\frac{\partial}{\partial t}\phi=u-a
\end{equation}
with initial value $\phi(0)=0$. Obviously $\phi(t)$ is a relative K\"ahler potential of $g(t)$ in the sense that \cite{Zhang}
\begin{equation}\label{e26}
g_{i\bar{j}}(t)=g_{i\bar{j}}(0)+\partial_i\partial_{\bar{j}}\phi(t).
\end{equation}

\begin{lemma}\label{l24}
Given $n$ and $\Lambda$ as in condition {\rm(\ref{C1})}, there
exists an increasing function
$\Phi:\mathbb{R}^+\rightarrow\mathbb{R}^+$ such that the following estimate holds under the K\"ahler-Ricci flow,
\begin{equation}\label{e27}
n+\triangle_0\phi(t)\leq\Phi(\alpha(t))
\end{equation}
where $\triangle_0$ denotes the Laplace of the initial metric $g_0$ and $\alpha(t)$ is defined by
\begin{equation}\nonumber
\alpha(t)=\sup\{\|\phi(s)\|_{C^0}+\|u(0)\|_{C^0}+\|u(s)\|_{C^0}|s\in[0,t]\}.
\end{equation}
\end{lemma}
\begin{proof}
Set $F(t)=\log\frac{\det g(t)}{\det g(0)}$. It is an easy check that
$$F(t)=u(t)-u(0)-\phi(t)-c(t),$$ for a family of constants $c(t)$ determined by
$$e^{c(t)}=\frac{1}{\V}\int_Me^{-\phi(t)-u(0)}dv_{g_0}.$$ As in the proof of Theorem 4.21 in \cite{ChLiWa}, combining Yau's estimate
\begin{eqnarray}
\triangle\big(e^{-k\phi}(n+\triangle_0\phi)\big)&\geq& e^{-k\phi}(\triangle_0F-s(g_0))-kne^{-k\phi}(n+\triangle_0\phi)\nonumber\\
&&+(k+\inf_{i,j}R_{i\bar{i}j\bar{j}}(g_0))e^{-k\phi-\frac{F}{n-1}}(n+\triangle_0\phi)^{\frac{n}{n-1}}\nonumber
\end{eqnarray}
and the evolution
\begin{eqnarray}\nonumber
\frac{\partial}{\partial t}\big(e^{-k\phi}(n+\triangle_0\phi)\big)=-k(u-a)e^{-k\phi}(n+\triangle_0\phi)+e^{-k\phi}\triangle_0(F+u(0)+\phi)
\end{eqnarray}
we get
\begin{eqnarray}
(\triangle-\frac{\partial}{\partial t})\big(e^{-k\phi}(n+\triangle_0\phi)\big)&\geq& e^{-k\phi}(n-\triangle_0u(0)-s(g_0))\nonumber\\
&&+(k(u-a)-kn-1)e^{-k\phi}(n+\triangle_0\phi)\nonumber\\
&&+(k+\inf_{i, j}R_{i\bar{i}j\bar{j}}(g_0))e^{-k\phi-\frac{F}{n-1}}(n+\triangle_0\phi)^{\frac{n}{n-1}}\nonumber\\
&=&(k(u-a)-kn-1)e^{-k\phi}(n+\triangle_0\phi)\nonumber\\
&&+(k+\inf_{i,j}R_{i\bar{i}j\bar{j}}(g_0))e^{-k\phi-\frac{F}{n-1}}(n+\triangle_0\phi)^{\frac{n}{n-1}}\nonumber.
\end{eqnarray}
Choose $k=\Lambda+1$. Fix any time interval $[0,t]$, suppose $e^{-k\phi}(n+\triangle_0\phi)$ achieve its maximal point at $(x_0,t_0)$, then
\begin{equation}
e^{-k\phi(x_0,t_0)}(n+\triangle_0\phi(x_0,t_0))\leq(1+kn-k(u(x_0,t_0)-a))^{n-1}e^{F(x_0,t_0)-k\phi(x_0,t_0)}.\nonumber
\end{equation}
Thus,
$$\sup_{[0,t]}(n+\triangle_0\phi)\leq\big(1+kn+k\|u-a\|_{C^0}(t_0)\big)^{n-1}e^{\sup(F(t_0))+k\osc(\phi(t_0))}.$$
This gives the desired bound of $(n+\triangle_0\phi)$ in view of the definition of $F$.
\end{proof}

By definition we have $\det g(t)=e^F\det g(0)$, consequently
$$\Phi^{1-n}e^F g_0\leq g(t)\leq\Phi g_0,$$
where $\Phi=\Phi(\alpha(t))$ for $\Phi$ and $\alpha(t)$ defined as in above lemma. After a modification to $\Phi$, we may assume that
\begin{equation}\label{e28}
\Phi^{-1}g_0\leq g(t)\leq\Phi g_0
\end{equation}
for $\Phi=\Phi(\alpha(t))$, an increasing function depending on $n$ and $\Lambda$ in Condition (\ref{C1}).

\begin{remark}
This is the only place where we essentially make use of the lower bound of the bisectional curvature. One might relax the condition on bisectional curvature to a condition of orthogonal holomorphic bisectional curvature.
\end{remark}

\subsection{Proof of Theorem \ref{main theorem}}

For simplicity, denote by $C_i,i=1,2,\cdots,$ a family of constants
depending on $n$, $\delta$, $\Lambda$ and $D$. The key step in
the proof of our Theorem \ref{main theorem} is the following:

\begin{claim}\label{c20}
Let $g(t)$ be the K\"ahler-Ricci flow with initial metric $g(0)$ satisfying {\rm(\ref{C1})} as above. There exists a constant $\Phi_0$ depending on $\delta$, $\Lambda$ and $D$ such that the following holds: for any $\epsilon$ small enough, one can find $L$ such that
\begin{equation}\label{e29}
\|\nabla u(t)\|_{C^0}< L\epsilon,\hspace{0.3cm}\Phi_0^{-1}g_0< g(t)<\Phi_0 g_0,\hspace{0.3cm}\lambda(g(t))> 1+\frac{\delta}{2}
\end{equation}
holds for all time $t\in[0,\infty)$.
\end{claim}

Let $L$ and $\Phi_0$ be given sufficiently large constants in a
prior, but $\Phi_0L\epsilon$ is small enough (see below for the
reason). Let $T$ be the largest time such that (\ref{e29}) holds on
the time interval $[0,T)$. Then, according to the preliminary lemmas
in \S 2.1, the K\"ahler-Ricci flow $g(t)$ behaves very well on
$[0,T)$. The aim is to show that $T=\infty$ for suitably chosen
$\epsilon$, $L$ and $\Phi_0$.

We argue by contradiction: suppose $T<\infty$ from now on, then we need to derive a contradiction for suitably chosen $\epsilon$, $L$ and $\Phi_0$.

First of all, by metric equivalence,
\begin{equation}\label{e210}
\diam(M,g(t))\leq\Phi_0^{1/2}D,\hspace{0.3cm}\forall t\in[0,T].
\end{equation}
Consequently
\begin{equation}\label{e211}
\|u-a\|_{C^0}\leq\Phi_0^{1/2}DL\epsilon,\hspace{0.3cm}\forall t\in[0,T].
\end{equation}
In particular,
\begin{equation}\label{e212}
\|u\|_{C^0}\leq 2\Phi_0^{1/2}DL\epsilon,\hspace{0.3cm}\forall t\in[0,T].
\end{equation}
In view of Lemma \ref{l22}, the first assumption we should make is
\begin{equation}\label{A1}
\Phi_0^{1/2}DL\epsilon\leq\min(\frac{1}{4},\frac{\delta_0^{'}}{8})
\end{equation}
where $\delta_0^{'}=\delta^{'}(\frac{\delta}{2},1)$ for
$\delta_0^{'}$ defined in (\ref{e20}), which is a constant depending
\textit{only} on $\delta$ in condition (\ref{C1}). Notice that when
$\Phi_0^{1/2}DL\epsilon\leq\frac{1}{4}$, we have by (\ref{e211}) that
$\osc(u)\leq2\Phi_0^{1/2}DL\epsilon\leq 1$ and thus
$\delta_0^{'}=\delta^{'}(\frac{\delta}{2},1)$ works for Lemma
\ref{l21}. Then Lemma \ref{l22} gives the exponential decay of
$Y(t)$:
\begin{equation}\label{e213}
Y(t)\leq e^{-\delta_0^{'}t}Y(0),\hspace{0.3cm}\forall t\in[0,T].
\end{equation}

To derive the $C^0$ estimate of $(u-a)$, we want to use Lemma
\ref{l23}. For this purpose, we recall that there exists $\kappa_0$
depending only on $\Lambda$ and $D$ such that $g_0$ is
$\kappa_0$-non-collapsed on scales $\leq 1$. As $g(t)$ is
$\Phi_0$-equivalent to $g_0$, we know that for any metric ball
$B_{g(t)}(r)$ of radius $r\leq\Phi_0^{-1/2}$ at time $t\leq T$, the
volume ratio has the bound
\begin{equation}\label{e214}
\frac{\Vol(B_{g(t)}(r))}{\V}\geq\kappa_0\Phi_0^{-2n}r^{-2n},\hspace{0.3cm}\forall t\in[0,T].
\end{equation}
Applying Lemma \ref{l23} for $\kappa=\kappa_0\Phi_0^{-2n}$ and $\rho=\Phi_0^{-1/2}$ we obtain
\begin{equation}\label{e215}
\|u-a\|_{C^0}\leq C_1\Phi_0^{\frac{n}{n+1}}Y^{\frac{1}{2n+2}},\hspace{0.3cm}\forall t\in[0,T],
\end{equation}
under the assumption that (to satisfy $\frac{\|u-a\|_{C^0}}{2}\leq\Phi_0^{-1/2}$), by (\ref{e211}),
\begin{equation}\label{A4}
\Phi_0DL\epsilon\leq 1.
\end{equation}
Combining with (\ref{e213}), (\ref{e215}) gives
\begin{equation}\label{e216}
\|u-a\|_{C^0}\leq C_1\Phi_0^{\frac{n}{n+1}}e^{-\frac{\delta^{'}_0t}{2n+2}}\epsilon^{\frac{1}{n+1}},\hspace{0.3cm}\forall t\in[0,T]
\end{equation}
since by weighted Poincar\'e inequality (\ref{RPP}),
$$Y(0)=\frac{1}{\V}\int_M(u-a)^2e^{-u}dv\leq\frac{1}{\V}\int_M|\nabla u|^2e^{-u}dv\leq\epsilon^2.$$

We first claim that the metric relation in (\ref{e29}) can be improved on time interval $[0,T]$. Define a function $\epsilon_0:\mathbb{R}^+\rightarrow\mathbb{R}^+$ via
\begin{equation}\label{A2}
\epsilon_0(A)= A^{-n}.
\end{equation}

\begin{claim}\label{c22}
There exist $\tilde{\Phi}$ depending on $\delta,\Lambda$ and $D$ with the following property. If $\Phi_0\geq 2\tilde{\Phi}$ and $\epsilon\leq\epsilon_0(\Phi_0)$ in our assumption {\rm(\ref{e29})}, then we have the estimate
\begin{equation}\label{e217}
\tilde{\Phi}^{-1}g_0\leq
g(t)\leq\tilde{\Phi}g_0,\hspace{0.3cm}\forall t\in[0,T].
\end{equation}
In particular,
\begin{equation}\label{e217}
2\Phi_0^{-1}g_0\leq
g(t)\leq\frac{1}{2}\Phi_0g_0,\hspace{0.3cm}\forall t\in[0,T].
\end{equation}
\end{claim}
\begin{proof}
In view of (\ref{e216}), if $\epsilon\leq\epsilon_0(\Phi_0)$, then $\|u-a\|_{C^0}\leq C_1e^{-\frac{\delta^{'}_0t}{2n+2}}$ on $[0,T]$. Thus by definition
$$\|\phi(t)\|_{C^0}\leq\int_0^t\|u-a\|_{C^0}ds\leq C_2,\hspace{0.3cm}\forall t\in[0,T].$$
By (\ref{e28}), this implies in particular that
$\tilde{\Phi}^{-1} g_0\leq g(t)\leq\tilde{\Phi} g_0$ on $[0,T]$ for
some constant $\tilde{\Phi}$ depending only on $n,\delta,\Lambda$
and $D$. The claim follows.
\end{proof}

\begin{claim}\label{c21}
Given $\Phi_0\geq 2\tilde{\Phi}$ sufficiently large. For any
$\epsilon\leq\epsilon_0(\Phi_0)$, set
\begin{equation}\label{L}
L=\Phi_0\epsilon^{\frac{1}{n+1}-1}
\end{equation}
in our assumption {\rm(\ref{e29})}. Then $g(t)$ satisfies
\begin{equation}\label{e218.5}
\|\nabla u\|_{C^0}(t)\leq
\frac{L}{2}\epsilon,\hspace{0.3cm}\forall t\in[0,T].
\end{equation}
\end{claim}
\begin{proof}
Recall the evolution of $|\nabla u|^2$ under the K\"ahler-Ricci flow \cite{SeTi}:
\begin{equation}\label{e36.5}
\frac{\partial}{\partial t}|\nabla u|^2=\triangle|\nabla u|^2-|\nabla\nabla u|^2-|\nabla\bar{\nabla}u|^2+|\nabla u|^2.
\end{equation}
From this we can bound $\|\nabla u\|_{C^0}$ on the short time interval $[0,2]$ as follows
\begin{equation}\label{e219}
\|\nabla u\|_{C^0}(t)\leq e\|\nabla u\|_{C^0}(0)\leq e\epsilon,\hspace{0.3cm}\forall t\in[0,2].
\end{equation}
When $t\geq 2$, we apply \cite[Page 234]{Ba} to derive
$$\|\nabla u\|_{C^0}(t)\leq e^2\|u-a\|_{C^0}(t-1).$$
Together with (\ref{e216}) this gives the
following bound:
\begin{eqnarray}
\|\nabla u\|_{C^0}(t)\leq C_1\Phi_0^{\frac{n}{n+1}}e^{-\frac{\delta^{'}_0(t-1)}{2n+2}}\epsilon^{\frac{1}{n+1}}
\leq C_3\Phi_0^{\frac{n}{n+1}}e^{-\frac{\delta^{'}_0t}{2n+2}}\epsilon^{\frac{1}{n+1}},\hspace{0.3cm}\forall t\in[2,T]\label{e220}.
\end{eqnarray}
Hence, under the assumption (\ref{L}), the estumate (\ref{e218.5})
is fulfilled whenever $\Phi_0$ is chosen large enough such that
\begin{equation}\label{Phi0}
\Phi_0^{\frac{1}{n+1}}\geq 2C_3
\end{equation}
for $C_3$ in (\ref{e220}). The proof is complete.
\end{proof}

From now on we assume $\Phi_0$ satisfies (\ref{Phi0}) and
$\Phi_0\geq2\tilde{\Phi}$ as in Claims \ref{c21} and \ref{c22}
respectively, where $\tilde{\Phi}$ is a constant depending only on $n,\delta,\Lambda$ and $D$. Let $\epsilon_0$ as in (\ref{A2}). Then choose $L$ according to the relation (\ref{L}) for
any given $\epsilon\leq\epsilon_0$. After these settings, conditions
(\ref{A1}) and (\ref{A4}) are equivalent to:
\begin{equation}\label{A1'}
\Phi_0^{\frac{3}{2}}D\epsilon^{\frac{1}{n+1}}\leq\min(\frac{1}{4},\frac{\delta^{'}_0}{8})\hspace{0.3cm}\mbox{and }\Phi_0^2D\epsilon^{\frac{1}{n+1}}\leq 1.
\end{equation}
Then $g(t)$ is $\tilde{\Phi}$-equivalent to $g(0)$ for any $t\in[0,T]$. Furthermore, putting (\ref{e219}) and (\ref{e220}) together we have the estimate
\begin{equation}\label{e221}
\|\nabla u\|_{C^0}(t)\leq C_4\Phi_0^{\frac{n}{n+1}}e^{-\frac{\delta^{'}_0t}{2n+2}}\epsilon^{\frac{1}{n+1}},\hspace{0.3cm}\forall t\in[0,T].
\end{equation}

We finally check the condition on eigenvalue $\lambda$ on the time interval $[0,T]$. Let us write $\lambda(t)=\lambda(g(t))$ for brevity.

\begin{claim}\label{c23}
Given $\Phi_0$ as above, there exists a constant $\epsilon_0^{'}=\epsilon_0^{'}(n,\Phi_0,\delta,\Lambda,D)$ such that if $\epsilon\leq\epsilon_0^{'}$ in our assumption, then we have the estimate on eigenvalue $\lambda(t)$:
\begin{equation}\label{e222}
\lambda(t)\geq 1+\frac{2}{3}\delta,\hspace{0.3cm}\forall t\in[0,T].
\end{equation}
\end{claim}
\begin{proof}
We may assume that $\lambda\leq 2$ on the time interval $[0,T]$. Let $\psi(t)$ be the corresponding eigenfunction of $\lambda(t)$ normalized by $\frac{1}{\V}\int_M|\psi|^2e^{-u}dv=1$. Then by definition, $-\triangle\psi+\nabla_iu\nabla_{\bar{i}}\psi=\lambda\psi$. Notice that by taking the time derivation to the normalization condition we get
$$\int_M(\dot{\psi}\bar{\psi}+\psi\bar{\dot{\psi}})e^{-u}dv=\int_M|\psi|^2(u-a)e^{-u}dv.$$ Then by an easy calculation,
\begin{eqnarray}
\frac{d}{dt}\lambda&=&\frac{d}{dt}\frac{1}{\V}\int_M|\bar{\nabla}\psi|^2e^{-u}dv\nonumber\\
&=&\frac{1}{\V}\int_M\big[-\nabla_i\nabla_{\bar{j}}u\nabla_{\bar{i}}\psi\nabla_j\bar{\psi}+\nabla_{\bar{i}}\dot{\psi}\nabla_i\bar{\psi}+
\nabla_{\bar{i}}\psi\nabla_i\bar{\dot{\psi}}-|\bar{\nabla}\psi|^2(u-a)\big]e^{-u}dv\nonumber\\
&=&\frac{1}{\V}\int_M\big[-\nabla_i\nabla_{\bar{j}}u\nabla_{\bar{i}}\psi\nabla_j\bar{\psi}+\lambda(\bar{\psi}\dot{\psi}+\psi\bar{\dot{\psi}})
-|\bar{\nabla}\psi|^2(u-a)\big]e^{-u}dv\nonumber\\
&=&\frac{1}{\V}\int_M\big[-\nabla_i\nabla_{\bar{j}}u\nabla_{\bar{i}}\psi\nabla_j\bar{\psi}+\lambda|\psi|^2(u-a)
-|\bar{\nabla}\psi|^2(u-a)\big]e^{-u}dv\nonumber\\
&\geq&-4\|u-a\|_{C^0}-\frac{1}{\V}\int_M\nabla_i\nabla_{\bar{j}}u\nabla_{\bar{i}}\psi\nabla_j\bar{\psi} e^{-u}dv.\nonumber
\end{eqnarray}
Thus, for any $t\in[0,T]$,
\begin{equation}\label{e223}
\lambda(t)\geq\lambda(0)-4\int_0^T\|u-a\|_{C^0}(s)ds-\frac{1}{\V}\int_0^T\int_M|\nabla_i\nabla_{\bar{j}}u\nabla_{\bar{i}}\psi\nabla_j\bar{\psi}| e^{-u}dvds.
\end{equation}
For our purpose, we need to estimate the integrals on the right hand side by universal constants independent of $T$. The first integral can be estimated by (\ref{e216}):
\begin{eqnarray}\label{e224}
\int_0^T\|u-a\|_{C^0}(s)ds\leq C_1\Phi_0^{\frac{n}{n+1}}\epsilon^{\frac{1}{n+1}}\int_0^Te^{-\frac{\delta^{'}_0s}{2n+2}}ds\leq C_5\Phi_0^{\frac{n}{n+1}}\epsilon^{\frac{1}{n+1}};
\end{eqnarray}
while the second integral can be estimated as follows: denote $\Psi=\sup_{[0,T]}\|\bar{\nabla}\psi\|_{C^0}$, then by Schwarz inequality,
\begin{eqnarray}
&&\frac{1}{\V}\int_0^T\int_M|\nabla_i\nabla_{\bar{j}}u\nabla_{\bar{i}}\psi\nabla_j\bar{\psi}| e^{-u}dvds
\leq\Psi^2\frac{1}{\V}\int_0^T\int_M|\nabla\bar{\nabla}u|e^{-u}dvds\nonumber\\
&\leq&\Psi^2\big(\sum_{i=0}^{[T]-1}\frac{1}{\V}\int_{i}^{i+1}\int_M|\nabla\bar{\nabla}u|e^{-u}dvds
+\int_{[T]}^T\frac{1}{\V}\int_{i}^{i+1}\int_M|\nabla\bar{\nabla}u|e^{-u}dvds\big)\nonumber\\
&\leq&\Psi^2\sum_{i=0}^{[T]-1}\big(\frac{1}{\V}\int_{i}^{i+1}\int_M|\nabla\bar{\nabla}u|^2e^{-u}dvds\big)^{1/2}
+\Psi^2\big(\frac{1}{\V}\int_{[T]}^T\int_M|\nabla\bar{\nabla}u|^2e^{-u}dvds\big)^{1/2}\label{e225}.
\end{eqnarray}
Thus, to estimate the second integral of (\ref{e223}), we need to consider two ingredients: (i) the uniform $L^\infty$ bound of $\bar{\nabla}\psi$ on the time interval $[0,T]$ (in terms of the Sobolev constant, as we will see later) and (ii) the $L^2$ integral bound of $|\nabla_i\nabla_{\bar{j}}u|$ on the space time. We shall deal with these two estimates one by one.

Firstly, we consider the $L^\infty$ uniform bound of $\bar{\nabla}\psi$ on $M\times[0,T]$. It is exactly the gradient estimate of eigenfunctions of $L$. To this aim, we should know the uniform bound of the Sobolev constant in a prior. By a well-known result of Croke \cite{Cr}, the isoperimetric constant, and thus the Sobolev constant \cite{ScYa}, of $g_0$ admits a universal bound by a constant depending only on $\Lambda$ and $D$. Since $g(t)$ is $\tilde{\Phi}$-equivalent to $g_0$ for any $t\in[0,T]$, there exists $C_6=C_6(n,\Lambda,D,\tilde{\Phi})=C_6(n,\delta,\Lambda,D)$ such that the Sobolev constant of $g(t)$, say $C_s(g(t))$, is uniformly bounded from above by $C_6$ for all $t\in[0,T]$. Consequently, by the gradient estimate, cf. Theorem \ref{gradient estimate}, there is a uniform constant $C_7$ such that
\begin{equation}\label{e226}
\Psi=\sup_{[0,T]}\|\bar{\nabla}\psi\|_{C^0}\leq C_7.
\end{equation}

Secondly, we estimate the integrals of $|\nabla_i\nabla_{\bar{j}}u|^2$ through studying the evolution of $Z=\frac{1}{\V}\int_M|\nabla u|^2e^{-u}dv$. Indeed, an easy calculation shows
\begin{eqnarray}
\frac{d}{dt}Z&=&\frac{1}{\V}\int_M\big[\triangle|\nabla u|^2-|\nabla\nabla u|^2-|\nabla\bar{\nabla}u|^2+|\nabla u|^2(1-u+a)\big]e^{-u}dv\nonumber\\
&=&\frac{1}{\V}\int_M\big[-\triangle u|\nabla u|^2+|\nabla u|^4-|\nabla\nabla u|^2-|\nabla\bar{\nabla}u|^2+|\nabla u|^2(1-u+a)\big]e^{-u}dv\nonumber\\
&\leq&\frac{1}{\V}\int_M\big[\frac{1}{2n}(\triangle u)^2-|\nabla\bar{\nabla}u|^2+(\frac{n}{2}+1)|\nabla u|^4+|\nabla u|^2(1-u+a)\big]e^{-u}dv\nonumber\\
&\leq&2Z-\frac{1}{\V}\int_M\frac{1}{2}|\nabla\bar{\nabla}u|^2e^{-u}dv,\nonumber
\end{eqnarray}
where we used $(\triangle u)^2\leq n|\nabla\bar{\nabla}u|^2$ and that $\|u-a\|_{C^0}$ and $\|\nabla u\|_{C^0}$ are small according to our assumption on the time interval $[0,T]$; see (\ref{e29}), (\ref{e211}) and assumption (\ref{A1}).
Now observe that $Z(t)$ can be estimated easily by $\|\nabla u\|_{C^0}$:
$$Z(t)=\frac{1}{\V}\int_M|\nabla u|^2e^{-u}dv\leq\|\nabla u\|_{C^0}^2(t).$$
Together with (\ref{e221}) we see the integration on any interval $[t,t+c]\subset[0,T]$, where $c\leq 1$, can be bounded as follows:
\begin{eqnarray}\label{e222}
\frac{1}{\V}\int_{t}^{t+c}\int_M|\nabla\bar{\nabla}u|^2e^{-u}dv&\leq& 2(Z(t)-Z(t+c))+4\int_{t}^{t+c}W(s)ds\nonumber\\
&\leq&2\|\nabla u\|_{C^0}^2(t)+4\int_t^{t+c}\|\nabla u\|_{C^0}^2(s)ds\nonumber\\
&\leq&2C_4^2\Phi_0^{\frac{2n}{n+1}}e^{-\frac{\delta^{'}_0t}{n+1}}\epsilon^{\frac{2}{n+1}}+4C_4^2\Phi_0^{\frac{2n}{n+1}}\epsilon^{\frac{2}{n+1}}
\int_t^{t+c}e^{-\frac{\delta_0^{'}s}{n+1}}ds\nonumber\\
&\leq&C_8\Phi_0^{\frac{2n}{n+1}}e^{-\frac{\delta^{'}_0t}{n+1}}\epsilon^{\frac{2}{n+1}}\label{e227}.
\end{eqnarray}

Now, substitute (\ref{e226}) and (\ref{e227}) into (\ref{e225}) we get the estimate of the second integral
\begin{eqnarray}\nonumber
\frac{1}{\V}\int_0^T\int_M|\nabla_i\nabla_{\bar{j}}u\nabla_{\bar{i}}\psi\nabla_j\bar{\psi}| e^{-u}dvds
\leq\Psi^2\sum_{i=0}^{[T]}C_8^{1/2}\Phi_0^{\frac{n}{n+1}}e^{-\frac{\delta^{'}_0i}{2n+2}}\epsilon^{\frac{1}{n+1}}\leq C_9\Phi_0^{\frac{n}{n+1}}\epsilon^{\frac{1}{n+1}}.
\end{eqnarray}
By (\ref{e223}), this estimate, together with (\ref{e224}), finally gives
\begin{eqnarray}
\lambda(t)&\geq&\lambda(0)-\int_0^T4\|u-a\|_{C^0}ds-\frac{1}{\V}\int_0^T\int_M|\nabla_i\nabla_{\bar{j}}u\nabla_{\bar{i}}\psi\nabla_j\bar{\psi}| e^{-u}dvds\nonumber\\
&\geq&\lambda(0)-4C_5\Phi_0^{\frac{n}{n+1}}\epsilon^{\frac{1}{n+1}}-C_9\Phi_0^{\frac{n}{n+1}}\epsilon^{\frac{1}{n+1}}\nonumber\\
&\geq&1+\delta-C_{10}\Phi_0^{\frac{n}{n+1}}\epsilon^{\frac{1}{n+1}}\nonumber.
\end{eqnarray}
Therefore, to make sure $\lambda(t)$ is bounded from below by $1+\frac{2}{3}\delta$, it only needs to add the following assumption on $\Phi_0$ and $\epsilon$:
\begin{equation}\label{A3}
C_{10}\Phi_0^{\frac{n}{n+1}}\epsilon^{\frac{1}{n+1}}\leq\frac{1}{3}\delta.
\end{equation}
The proof of Claim \ref{c23} is completed by setting $\epsilon^{'}_0$ such that the equality in (\ref{A3}) holds.
\end{proof}

\begin{proof}[Proof of Claim \ref{c20}]
First of all, choose any $\Phi_0$ satisfying (\ref{Phi0}) in Claim \ref{c21} and
$\Phi_0\geq2\tilde{\Phi}$ where $\tilde{\Phi}$ is the constant in
Claim \ref{c22}, which depends only on the constants $n,\delta,\Lambda$ and
$D$. Let $\epsilon_0^{''}$ be the largest constant such that
(\ref{A1'}) and (\ref{A3}) hold and
$\epsilon_0^{''}\leq\epsilon_0(\Phi_0)$ as in Claim \ref{c22}. Finally for any
$\epsilon\leq\epsilon_0^{''}$, choose $L$ as in (\ref{L}). Then, as
showed in above three claims, if the condition (\ref{e29}) holds on
any time interval $[0,T]$, then the following more strict condition
$$2\Phi_0^{-1}g(0)\leq g(t)\leq\frac{1}{2}\Phi_0g(0),\hspace{0.3cm}\|\nabla u\|_{C^0}\leq\frac{L}{2}\epsilon,\hspace{0.3cm}\lambda\geq 1+\frac{2}{3}\delta$$
still holds on the same interval $[0,T]$. Thus, the maximal time $T$ can only happen to be $\infty$. This proves the claim.
\end{proof}

\begin{proof}[Proof of Theorem \ref{main theorem}]
Fix a sufficiently large constant $\Phi_0$, choose $\epsilon$ and
$L$ such that Claim \ref{c20} holds and then apply Theorem 1.5 of
\cite{Zhang}.
\end{proof}

\section{Proof of Theorem \ref{t2}}

Let $(M,g_0,J)$ be a compact K\"ahler manifold of dimension $n$ with
$c_1(M;J)>0$. Suppose $g_0$ satisfies
\begin{equation}\label{e302}
|Rm(g_0)|\leq\Lambda,\hspace{0.2cm}\lambda(g_0)\geq
1+\delta,\hspace{0.2cm}\frac{1}{\V}\int_M(s-n)^2dv\leq\epsilon,
\end{equation}
for fixed $\delta,\Lambda$ and certain $\epsilon$ small enough,
where $\V$ is the volume of $g_0$ while $\lambda(g_0)$ is the
``second" eigenvalue of $g_0$ as defined in the Introduction.

Consider the K\"ahler-Ricci flow $g(t)$ starting from $g_0$ on $M$.
By the evolution of curvature tensor under the Ricci flow, there
exists $T_0=T_0(\Lambda)\leq 1$ such that
\begin{equation}\label{e31}
|Rm|(t)\leq 2\Lambda,\hspace{0.3cm}\forall t\in[0,T_0].
\end{equation}
We are going to show that the condition (\ref{Z1}) is fulfilled at
some time $t_0\leq T_0$ for certain constants
$\delta^{'},\Lambda^{'},D^{'}$ depending on $\delta,\Lambda$ and
sufficiently small $\epsilon^{'}$ depending on $\epsilon$. By
Theorem \ref{main theorem} it completes the proof of Theorem
\ref{t2}.

For simplicity, let us denote $\tilde{C}_i,i=1,2,\cdots,$ a family
of positive constants depending on $n,\delta$ and $\Lambda$.

We start with recalling some geometrical estimates on the time
interval $[0,T_0]$. First of all, a similar argument as in \cite[Lemma 3.5]{ChTi2} gives a bound of the diameter of the initial metric
$g_0$ by a constant $\tilde{C}_1$. Then, as each metric $g(t)$ is
$e^{2\Lambda}$-equivalent to $g_0$ (from the evolution of
K\"ahler-Ricci flow (\ref{KRF})), we get a uniform bound on the
diameter of the manifold:
\begin{equation}\label{e32}
\diam(M,g(t))\leq \tilde{C}_2,\hspace{0.3cm}\forall t\in[0,T_0].
\end{equation}
Then, by virtue of Li and Yau's eigenvalue estimate \cite{LiYa,
ScYa}, the first eigenvalue of Laplace with respect to each metric
$g(t)$ admits a uniform lower bound
\begin{equation}\label{e33}
\lambda_1(g(t))\geq\tilde{C}_3^{-1}.
\end{equation}
Or, in other words, the Poincar\'e inequality uniformly holds on the
time interval $[0,T_0]$. Furthermore, since the volume of the
K\"ahler manifold $(M,g(t))$ has a uniform lower bound depending
only on the dimension $n$, by a result of Bando and Mabuchi
\cite{BaMa}, the Green function of $(M,g(t))$, say $G_t(x,y)$ with
$\inf G_t=0$, has a uniform upper bound in the sense of integral
\begin{equation}\label{e34}
\int_MG_t(x,y)dv_{g(t)}(y)\leq \tilde{C}_4.
\end{equation}
Finally, by Croke's estimate \cite{Cr} on the isoperimetric constant
and the well-known equivalence of the isoperimetric inequality and
the Sobolev inequality \cite{ScYa}, the Sobolev constant of $(M,g(t))$, say
$C_{s,t}$, also admits a uniform bound:
\begin{equation}\label{e35}
C_{s,t}\leq\tilde{C}_5,\hspace{0.3cm}\forall t\in[0,T_0].
\end{equation}

The following argument can be divided into several lemmas.

\begin{lemma}
Let $u(t)$ be the Ricci potential of $g(t)$ determined by
{\rm(\ref{RP})} and {\rm(\ref{RPN})}. There is a uniform bound of
$u(t)$:
\begin{equation}\label{e38}
\|u(t)\|_{C^0}\leq \tilde{C}_6,\hspace{0.3cm}\forall t\in[0,T_0].
\end{equation}
\end{lemma}
\begin{proof}
Denote $(M,g)=(M,g(t))$ and $u=u(t)$ for $t\in[0,T_0]$. First of
all, the Jensen inequality gives
$$\frac{1}{\V}\int_M-udv\leq\log\big(\frac{1}{\V}\int_Me^{-u}dv\big)=0.$$
Then apply the Green formula, by (\ref{e34}), for any $x\in M$,
\begin{equation}\nonumber
u(x)=\frac{1}{\V}\int_Mu(y)dv(y)-\frac{1}{\V}\int_M\triangle u(y)G(x,y)dv(y)\geq-\tilde{C}_6
\end{equation}
since $-\triangle u=s-n\geq-2n^2\Lambda-n$ by (\ref{e31}). Now,
define $\varphi=u-\inf u+1$; to prove the $L^\infty$ bound of $u$,
it suffices to show a uniform upper bound of $\varphi$. This
follows from a classical iteration argument, cf. \cite{Ti1} for
example. Namely, one first obtains the $L^1$ bound of $\varphi$
by applying the Green formula at a maximal point of $\varphi$, say
$x_0$,
$$\frac{1}{\V}\int_M\varphi dv=1+\frac{1}{\V}\int_M\triangle\varphi(y)G(x_0,y)dv(y)\leq1+(2n^2\Lambda+n)\tilde{C}_4;$$
then apply the Poincar\'e inequality to get the $L^2$ bound of
$\varphi$:
\begin{eqnarray}
\lambda_1\int_M\varphi^2dv&\leq&\int_M|\nabla\varphi|^2dv+\big(\int_M\varphi dv\big)^2\nonumber\\
&=&-\int_M\varphi\triangle\varphi dv+\big(\int_M\varphi dv\big)^2\nonumber\\
&\leq&\big(\int_M\varphi^2 dv\big)^{1/2}\big(\int_M(\triangle\varphi)^2 dv\big)^{1/2}+\big(\int_M\varphi dv\big)^2\nonumber\\
&\leq&\frac{\lambda_1}{2}\int_M\varphi^2 dv+\frac{1}{2\lambda_1}\int_M(s-n)^2dv+\big(\int_M\varphi dv\big)^2\nonumber;
\end{eqnarray}
finally adopting Moser's iteration argument to
\begin{equation}\nonumber
\triangle\varphi=\triangle u\geq -(2n^2\Lambda+n)
\end{equation}
one deduces an upper bound of $\varphi$ in terms of $n,\Lambda,
C_{s,t}$ and $\int_M\varphi^2 dv$. Combining with these three estimates finishes proof
of the lemma.
\end{proof}

\begin{lemma}
There exists $t_0\leq T_0$ depending on $n,\delta$ and $\Lambda$ such
that
\begin{equation}\label{e39}
\lambda(g(t))\geq 1+\frac{\delta}{2},\hspace{0.3cm} \forall t\leq t_0.
\end{equation}
\end{lemma}
\begin{proof}
Let $\psi(t)$ be associated eigenfunctions of $\lambda(g(t))$ which
are normalized such that $\frac{1}{\V}\int_M|\psi|^2e^{-u}dv=1$. By
the calculation in the proof of Claim \ref{c23} in \S 2 we have
\begin{eqnarray}
\frac{d}{dt}\lambda&=&\frac{1}{\V}\int_M\big[-\nabla_i\nabla_{\bar{j}}u\nabla_{\bar{i}}\psi\nabla_j\psi+\lambda\psi^2(u-a)-|\nabla\psi|^2(u-a)\big]e^{-u}dv\nonumber\\
&\geq&-2\lambda\|u-a\|_{C^0}-\lambda\|\nabla_i\nabla_{\bar{j}}u\|_{C^0}\nonumber\\
&\geq&-4\tilde{C}_6\lambda-(n\Lambda+n)\lambda\nonumber
\end{eqnarray}
where we used (\ref{e38}) and that $\|\nabla_i\nabla_{\bar{j}}u\|_{C^0}=\|g_{i\bar{j}}-R_{i\bar{j}}\|_{C^0}\leq(n\Lambda+n)$ in the last inequality. Since
$\lambda(g(0))\geq 1+\delta$, we can choose a small $t_0$ depending
on $n,\delta$ and $\Lambda$ such that (\ref{e39}) holds.
\end{proof}

\begin{lemma}
$\|\nabla u\|_{C^0}(t_0)$ can be estimated as follows:
\begin{equation}\label{e310}
\|\nabla u\|_{C^0}(t_0)\leq\tilde{C}_8\epsilon^{1/4}.
\end{equation}
\end{lemma}
\begin{proof}
We will use the Moser iteration to prove this lemma; the details are referred to \cite{ChTi2}. From
the evolution of $|\nabla u|^2$, (\ref{e36.5}), we get the inequality:
\begin{equation}\nonumber
\frac{\partial}{\partial t}|\nabla u|^2\leq\triangle|\nabla
u|^2+|\nabla u|^2.
\end{equation}
To run the iteration argument, we need to (i) estimate the integral
of $|\nabla u|^2$ on the space time $M\times[0,t_0]$ and (ii) get
the uniform bound of $|\nabla u|$ on the time interval
$[\frac{t_0}{2},t_0]$. We first deal with (i). At the initial time
$t=0$, as $u(0)$ is uniformly bounded (\ref{e38}), the estimate can
be deduced as follows
\begin{eqnarray}
\frac{1}{\V}\int_M|\nabla u|^2dv&\leq& e^{\tilde{C}_6}\frac{1}{\V}\int_M|\nabla u|^2e^{-u}dv=e^{\tilde{C}_6}\frac{1}{\V}\int_M\triangle ue^{-u}dv\nonumber\\
&\leq& e^{\tilde{C}_6}\big(\frac{1}{\V}\int_M(\triangle
u)^2e^{-u}dv\big)^{1/2}\leq
e^{\tilde{C}_6+\sqrt{\tilde{C}_6}}\epsilon^{1/2}.\nonumber
\end{eqnarray}
For a general time $t\in[0,t_0]$ where $t_0\leq 1$, the estimate of
$\int_M|\nabla u|^2dv$ follows from the evolution of $|\nabla u|^2$
under the K\"ahler-Ricci flow. Actually,
\begin{eqnarray}
\frac{d}{dt}\int_M|\nabla u|^2dv &=&\int_M\big(\triangle|\nabla
u|^2-|\nabla u|^2-|\nabla\bar{\nabla}u|^2+|\nabla u|^2+\triangle
u|\nabla
u|^2\big)dv\nonumber\\
&\leq&(1+n+2n^2\Lambda)\int_M|\nabla u|^2dv,\nonumber
\end{eqnarray}
where we used $|\triangle u|=|n-s|\leq n+2n^2\Lambda$ on $[0,T_0]$.
Thus,
\begin{equation}\label{e311}
\frac{1}{\V}\int_M|\nabla u|^2dv\leq
e^{1+n+2n^2\Lambda}e^{\tilde{C}_6+\sqrt{\tilde{C}_6}}\epsilon^{1/2},\hspace{0.3cm}\forall
t\in[0,t_0].
\end{equation}

We next deal with (ii). This actually follows from a smoothing argument due to Bando \cite{Ba}; see \cite[Lemma 1]{PSSW1}
also. Indeed, from the argument there, one knows that
\begin{equation}
\|\nabla u\|_{C^0}(t)\leq
\big(\frac{2}{t_0}\big)^{-1/2}e\|u\|_{C^0}(0)\leq\tilde{C}_9,\hspace{0.3cm}\forall
t\in[\frac{t_0}{2},t_0].\label{e312}
\end{equation}

Now, since the Poincar\'e inequality and Sobolev inequality hold
uniformly on the time interval $[0,t_0]$, an iteration argument
\cite{ChTi1, ChLi} shows
\begin{eqnarray}
\|\nabla
u\|_{C^0}^2(t_0)\leq\tilde{C}_{10}t_0^{-\frac{n+2}{4}}\big(\int_{\frac{t_0}{2}}^{t_0}\int_M|\nabla
u|^4dvdt\big)^{1/2}.\label{e313}
\end{eqnarray}
Combining with (\ref{e311}) and (\ref{e312}) this gives the desired
result.
\end{proof}

Once known these lemmas, the proof of Theorem \ref{t2} is clear.

\begin{proof}[Proof of Theorem \ref{t2}]
Putting the estimates (\ref{e31}), (\ref{e32}), (\ref{e39}) and
(\ref{e310}) together, we know $g(t_0)$ satisfies the following
condition:
$$|Rm|(t_0)\leq2\Lambda,\hspace{0.3cm}\diam(M,g(t_0))\leq\tilde{C}_1,\hspace{0.3cm}\lambda(g(t_0))\geq 1+\frac{\delta}{2},
\hspace{0.3cm}\|\nabla u\|_{C^0}(t_0)\leq\tilde{C}_8\epsilon^{1/4}$$
at some $t_0>0$ where $\tilde{C}_1$ and $\tilde{C}_7$ are constants depending only
on $n,\delta$ and $\Lambda$ in (\ref{Z2}). By Theorem \ref{main
theorem}, the K\"ahler-Ricci flow will converge to a
K\"ahler-Einstein metric whenever $\epsilon$ is small enough. This
completes the proof of the Theorem \ref{t2}.
\end{proof}

\section{Further remarks and questions on stability}

\begin{remark}
Let $(M,g_{KE},J_{KE})$ be a compact K\"ahler-Einstein manifold with positive scalar curvature and $(M,J)$ be an arbitrarily small deformed K\"ahler manifold with $c_1(M;J)>0$. If $h^0(J_{KE})=0$, then by the argument of Lemma 1.3 in \cite{Ti0}, $(M,J)$ admits a K\"ahler-Einstein metric. This can also be derived by \cite[Theorem 1.4]{ChLi} or our Theorem 1.8. In the case $\dim h^0(J_{KE})>0$, Chen and Li proved \cite{ChLi} $(M,J)$ admits a K\"ahler-Einstein metric if $J$ is pre-stable and has vanished Futaki invariant on $\pi c_1(M;J)$. Sz\'ekelyhidi proved in \cite[Theorem 2]{Sz} the stability of existence of constant scalar curvature metrics under the K-polystable conditions, which shows in particular the existence of K\"ahler-Einstein metric on $(M,J)$ if it is K-polystable. Notice that K-polystable is also ``almost" necessary for the existence of K\"ahelr-Einstein metric in $\pi c_1(M;J)$, see \cite{St, StSz}, thus Sz\'ekelyhidi's theorem almost solves the stability problem completely. Here, our Theorem \ref{t4} give another insight to this problem in view of Tian and Donaldson's counterexample.
\end{remark}

\begin{remark}
One natural question is how to extend our stability theorem to the case of shrinking K\"ahler-Ricci solitons. However, two problems arise in the later situation: (1) how to determine a suitable holomorphic vector field with respect to which the K\"ahler-Ricci soliton structure exists? (2) how to relate the K\"ahler-Ricci flow on the background K\"ahler-Ricci soliton and the K\"ahler-Ricci flow on the perturbed K\"ahler manifold on which the K\"ahler-Ricci soliton may exists? To the author's understanding, these may not be easy problems.
\end{remark}

The assumption of vanished Futaki invariant in the stability Theorem \ref{t3} and \ref{t4} is a closed condition on complex structures and so looks not good. Comparing Tian's argument \cite{Ti0}, we ask

\begin{question}
Can we remove the condition of vanished Futaki invariant in Theorem \ref{t3} and \ref{t4}? Or, in other words, does the Futaki invariant vanish automatically under other conditions?
\end{question}

As known to experts in the field of K\"ahler geometry, the solvability of the existence of K\"ahler-Einstein metric for a generic complex structure can hopefully be reduced to prove the closedness of the space of K\"ahler-Einstein metrics under certain stability condition, such as the K-stability or other stabilities in GIT; see Tian's celebrated solution to the existence of K\"ahler-Einstein metrics on complex surfaces \cite{Ti0} and Donaldson's note \cite{Do2} for higher dimensional case. According to this discussion, we ask a similar question in our setting.

\begin{question}
Is the space of K\"ahler-Einstein metrics closed under the ``relative stable" condition? More precisely, let $(M,g_\infty,J_\infty)$ be a compact K\"ahler manifold whose K\"ahler form lies in $\pi c_1(M;J_\infty)$ and $(g_i,J_i)$ be a sequence of K\"ahler-Einstein metrics on $M$ with positive scalar curvature such that $J_i\rightarrow J_\infty$ in the $C^2$ topology, if $J_\infty$ is ``relative stable" (and possibly the Futaki invariant vanishes on $\pi c_1(M;J_\infty)$), then does there exist a K\"ahler-Einstein metric on $(M,J_\infty)$?
\end{question}

The following question may have been well understood by experts, but since the lack of knowledge, the author still write it down here as a problem.

\begin{question}
Establish the relation between the ``relative stable" and other stability conditions from GIT.
\end{question}

\section{Appendix. Gradient estimate on K\"ahler manifolds with $c_1>0$}

In this section, we will use the Moser's iteration argument to prove
the gradient estimate to the eigenfunctions of the operator
$L=-\triangle+g^{i\bar{j}}\nabla_iu\nabla_{\bar{j}}$ on a K\"ahler
manifold with positive first Chern class. The approach is classical when considering the gradient estimate for eigenfunctions of Laplace on a manifold with Ricci curvature bounded from below.

\begin{theorem}\label{gradient estimate}
Let $(M,g,J)$ be a compact K\"ahler manifold of dimension $n$ with
$c_1(M;J)>0$ and $u$ be the Ricci potential determined by
{\rm(\ref{RP})} and {\rm(\ref{RPN})}. Suppose $\psi\in C^\infty(M)$
satisfies
\begin{equation}\label{eigenfunction}
-\triangle\psi+g^{i\bar{j}}\nabla_iu\nabla_{\bar{j}}\psi=\lambda\psi
\end{equation}
for some $\lambda>0$. Then we have
\begin{equation}\label{GA}
\sup|\bar{\nabla}\psi|\leq
C(n)e^{\frac{\max(u)}{2}}C_s^{n/2}(\|\nabla
u\|_{C^0}^2+\lambda)^{n/2}\lambda^{1/2}
\end{equation}
where $C(n)$ is a constant depending on $n$, while $C_s$ denotes the Sobolev constant of $g$ such that
\begin{equation}\label{Sobolev}
\big(\int_M|\nabla f|^{\frac{2n}{n-1}}dv\big)^{\frac{n-1}{n}}\leq
C_s\int_M\big(|\nabla f|^2+f^2\big)dv,\hspace{0.3cm}\forall f\in
C^\infty(M;\mathbb{R}).
\end{equation}
\end{theorem}
\begin{proof}
We start with the Bochner formula for $f=|\bar{\nabla}\psi|^2$:
\begin{eqnarray}\nonumber
\triangle f&=&|\bar{\nabla}\bar{\nabla}\psi|^2+|\nabla\bar{\nabla}\psi|^2+\nabla_j\triangle\bar{\psi}\nabla_{\bar{j}}\psi
+\nabla_{\bar{j}}\triangle\psi\nabla_j\bar{\psi}+R_{i\bar{j}}\nabla_{\bar{i}}\psi\nabla_j\bar{\psi}\nonumber\\
&=&|\bar{\nabla}\bar{\nabla}\psi|^2+|\nabla\bar{\nabla}\psi|^2+(1-2\lambda)f+\nabla_i\nabla_{\bar{j}}u\nabla_{\bar{i}}\psi\nabla_j\bar{\psi}\nonumber\\
&&+\nabla_i\nabla_j\bar{\psi}\nabla_{\bar{i}}u\nabla_{\bar{j}}\psi+\nabla_{\bar{i}}\nabla_{\bar{j}}\psi\nabla_iu\nabla_j\bar{\psi}\nonumber\\
&\geq&\frac{1}{2}|\bar{\nabla}\bar{\nabla}\psi|^2+|\nabla\bar{\nabla}\psi|^2+(1-2\lambda-|\nabla u|^2)f+\nabla_i\nabla_{\bar{j}}u\nabla_{\bar{i}}\psi\nabla_j\bar{\psi}.\nonumber
\end{eqnarray}
where we used the Schwarz inequality in the last inequality. Then, for any $q\geq 2$,
\begin{eqnarray}
&&-\frac{4(q-1)}{q^2}\int|\nabla f^{q/2}|^2dv=\int f^{q-1}\triangle fdv\nonumber\\
&\geq&\int f^{q-1}\big(\frac{1}{2}|\bar{\nabla}\bar{\nabla}\psi|^2+|\nabla\bar{\nabla}\psi|^2+(1-2\lambda-|\nabla u|^2)f+\nabla_i\nabla_{\bar{j}}u\nabla_{\bar{i}}\psi\nabla_j\bar{\psi}\big)dv.\label{Ap2}
\end{eqnarray}
The last term on the right hand side can be estimated as follows:
\begin{eqnarray}
&&\int f^{q-1}\nabla_i\nabla_{\bar{j}}u\nabla_{\bar{i}}\psi\nabla_j\bar{\psi}dv\nonumber\\
&=&
\int\big(-(q-1)f^{q-2}\nabla_if\nabla_{\bar{i}}\psi\nabla_{\bar{j}}u\nabla_{j}\bar{\psi}\nonumber\\
&&-f^{q-1}\triangle\psi\nabla_{\bar{j}}u\nabla_j\bar{\psi}-f^{q-1}\nabla_i\nabla_j\bar{\psi}\nabla_{\bar{i}}\psi\nabla_{\bar{j}}u
\big)dv\nonumber\\
&=&\int\big(-\frac{2(q-1)}{q}f^{\frac{q}{2}-1}\nabla_if^{q/2}\nabla_{\bar{i}}\psi\nabla_{\bar{j}}u\nabla_j\bar{\psi}\nonumber\\
&&-f^{q-1}\triangle\psi\nabla_{\bar{j}}u\nabla_j\bar{\psi}
-f^{q-1}\nabla_i\nabla_j\bar{\psi}\nabla_{\bar{i}}\psi\nabla_{\bar{j}}u\big)dv\nonumber\\
&\geq&-\frac{q-1}{q^2}\int|\nabla f^{q/2}|^2dv-(q-1)\|\nabla
u\|_{C^0}^2\int f^qdv\nonumber\\
&&-\frac{1}{n}\int f^{q-1}(\triangle\psi)^2dv-\frac{n}{4}\|\nabla u\|^2_{C^0}\int f^qdv\nonumber\\
&&-\frac{1}{2}\int
f^{q-1}|\bar{\nabla}\bar{\nabla}\psi|^2dv-\frac{1}{2}\|\nabla
u\|_{C^0}^2\int f^qdv\nonumber\\
&\geq&-\frac{q-1}{q^2}\int|\nabla
f^{q/2}|^2dv-(q+\frac{n}{4})\|\nabla u\|_{C^0}^2\int
f^qdv\nonumber\\
&&-\frac{1}{n}\int f^{q-1}(\triangle\psi)^2dv-\frac{1}{2}\int
f^{q-1}|\bar{\nabla}\bar{\nabla}\psi|^2dv\nonumber,
\end{eqnarray}
Substituting it into (\ref{Ap2}) we obtain
\begin{eqnarray}
\frac{3(q-1)}{q^2}\int|\nabla f^{q/2}|^2dv&\leq&[(q+C)\|\nabla
u\|_{C^0}^2+2\lambda-1]\int f^qdv\nonumber,
\end{eqnarray}
where $C=C(n)$. We mention that we have used $(\triangle\psi)^2\leq n|\nabla\bar{\nabla}\psi|^2$ here. Thus
$$\int|\nabla f^{q/2}|^2dv\leq \big(Cq^2(\|\nabla u\|_{C^0}^2+\lambda)-1\big)\int f^qdv$$
for a uniform constant $C$ depending only on $n$. By Sobolev
inequality (\ref{Sobolev}),
\begin{eqnarray}\nonumber
\big(\int f^{q\mu}dv\big)^{1/\mu}\leq C_sCq^2(\|\nabla
u\|_{C^0}^2+\lambda)\int f^qdv,\hspace{0.3cm}\forall q\geq 2,
\end{eqnarray}
where $\mu=\frac{n}{n-1}$. Define $q_k=2\mu^k$ for $k\geq 0$, then
$$\|f\|_{L^{q_{k+1}}}\leq[C_sC(\|\nabla u\|_{C^0}^2+\lambda)]^{1/q_k}q_k^{2/q_k}\|f\|_{L^{q_k}},\hspace{0.3cm}\forall k\geq 0.$$
By induction we have
$$\sup f\leq C(n)C_s^{\frac{n}{2}}(\|\nabla u\|_{C^0}^2+\lambda)^{\frac{n}{2}}\|f\|_{L^2}.$$
Since $\|f\|_{L^2}\leq(\sup f)^{1/2}\|f\|_{L^1}^{1/2}$, this implies
immediately
\begin{eqnarray}
\sup f&\leq& C(n)C_s^n(\|\nabla u\|_{C^0}^2+\lambda)^{n}\int fdv\nonumber\\
&=&C(n)C_s^n(\|\nabla
u\|_{C^0}^2+\lambda)^{n}\int|\bar{\nabla}\psi|^2dv\nonumber.
\end{eqnarray}
Notice that $\int|\bar{\nabla}\psi|^2e^{-u}dv=\lambda$, we finally
obtain
\begin{eqnarray}
\sup f&\leq&C(n)e^{\max(u)}C_s^n(\|\nabla u\|_{C^0}^2+\lambda)^{n}\int|\bar{\nabla}\psi|^2e^{-u}dv\nonumber\\
&=&C(n)e^{\max(u)}C_s^n(\|\nabla
u\|_{C^0}^2+\lambda)^{n}\lambda.\nonumber
\end{eqnarray}
The result now follows.
\end{proof}

\end{document}